\documentclass[letterpaper, 10 pt, conference]{ieeeconf}%
\usepackage{graphics}
\usepackage{epsfig}
\usepackage{cite}
\usepackage{breakurl}
\usepackage{amsmath}
\usepackage{amssymb}
\usepackage{amsfonts}
\usepackage{graphicx}%
\providecommand{\U}[1]{\protect\rule{.1in}{.1in}}
\IEEEoverridecommandlockouts
\newtheorem{theorem}{Theorem}[section]

\newtheorem{assumption}[theorem]{Assumption}

\allowdisplaybreaks[0]

\usepackage[normalem]{ulem}
\usepackage{color}
\usepackage{dsfont}
\usepackage[hyphens]{url}
\usepackage[hidelinks]{hyperref}
\hypersetup{breaklinks=true}
\begin{document}

\title{{\LARGE \textbf{Decentralized Optimal Control for Connected Automated Vehicles at Intersections Including Left and Right Turns}}}
\author{ Yue Zhang, Andreas A. Malikopoulos, Christos G. Cassandras
\thanks{This research was supported by US
Department of Energy's SMART Mobility Initiative. The work of Cassandras and
Zhang is supported in part by NSF under grants ECCS-1509084, CNS-1645681 and IIP-1430145, by AFOSR under grant FA9550-15-1-0471, and by a grant from the
MathWorks.} \thanks{Y. Zhang and C.G. Cassandras are with the Division of
Systems Engineering and Center for Information and Systems Engineering, Boston
University, Boston, MA 02215 USA (e-mail: joycez@bu.edu; cgc@bu.edu).}
\thanks{A.A. Malikopoulos is with the Department of Mechanical Engineering, University of Delaware, Newark, DE 19716 USA (email: andreas@udel.edu).} }
\maketitle

\begin{abstract}
In prior work, we addressed the problem of optimally controlling on line connected and automated vehicles crossing two adjacent intersections in an urban area to minimize fuel consumption while achieving maximal throughput without any explicit traffic signaling and without considering left and right turns. In this paper, we extend the solution of this problem  to account for left and right turns under hard safety constraints. Furthermore, we formulate and solve another optimization problem to minimize a measure of passenger discomfort while the vehicle turns at the intersection and we investigate the associated tradeoff between minimizing fuel consumption and passenger discomfort.
\end{abstract}

\thispagestyle{empty} \pagestyle{empty}


\section{Introduction}
Traffic light signaling is the prevailing method used to control the traffic flow through an
intersection.  Aside from the infrastructure cost and the need for
dynamically controlling green/red cycles, traffic light systems 
 can increase
the number of rear-end collisions at the intersection \cite{Ohio}. Serious
delays can occur during hours of heavy traffic if the light
cycle is not adjusted appropriately.  These challenges have motivated research efforts for new approaches capable of providing a smoother traffic flow and more
fuel-efficient driving while also improving safety. 

 Dresner and Stone \cite{Dresner2004} proposed the use of a centralized reservation scheme to control a single intersection of two roads with vehicles traveling with similar speed on a single direction on each road, i.e., no turns are allowed. Since then, several other efforts using reservation schemes have been reported in the literature \cite{Dresner2008,DeLaFortelle2010,Huang2012}. Increasing the throughput of an intersection is one desired goal which can be
achieved through the travel time optimization of all vehicles located within a
radius from the intersection. Several efforts have focused on minimizing
vehicle travel time under collision-avoidance constraints
\cite{Li2006,Yan2009,Zohdy2012, Zhu2015}. Lee and Park \cite{Lee2012} proposed
an approach based on minimizing the overlap in the position of
vehicles inside the intersection rather than their arrival times. 
A detailed discussion of the research in this area reported
in the literature to date can be found in \cite{Rios-Torres}.

Connected and automated vehicles (CAVs) provide the most intriguing and promising opportunity to reduce fuel consumption, greenhouse gas emissions, travel delays and to improve safety. In earlier work \cite{ZhangMalikopoulosCassandras2016}, we established a decentralized optimal control framework to address the problem of optimally controlling CAVs crossing two adjacent intersections in an urban area, without any explicit traffic signaling and without considering left and right turns, with the objective of minimizing fuel consumption while achieving maximal throughput using a First-In-First-Out queue to designate the order in which the CAVs cross the intersection. We also established the conditions under which feasible solutions exist and showed that they can be enforced through an appropriately designed \textit{Feasibility Enforcement Zone} (FEZ) that precedes the \textit{Control Zone} (CZ) in \cite{Zhang2016}.

To consider left and right turns within this framework, ``comfort" becomes of fundamental importance in addition to safety. In this paper, we extend the solution of the problem addressed in \cite{ZhangMalikopoulosCassandras2016} to account for left and right turns. Then, another optimization problem is formulated with the objective of minimizing a measure of passenger discomfort while the vehicle turns.
Furthermore, we investigate the associated tradeoff between minimizing fuel consumption and passenger discomfort inside the \textit{Merging Zone} (MZ). 

The problem of coordinating CAVs at intersections including left and right turns has been addressed in \cite{kim2014mpc} using an approach based on Model Predictive Controlto to achieve system-wide safety and liveness of intersection-crossing traffic. However, in our approach, the objective is to jointly minimize fuel consumption and passenger discomfort for CAVs crossing an intersection, while safety is a hard constraint.

The paper is organized as follows. In Section II, we review the model
in \cite{ZhangMalikopoulosCassandras2016} and its generalization in
\cite{Malikopoulos2016}. In Section III, we present the model and the analytical solution of the decentralized optimal control problem with the new collision-avoidance terminal conditions. In Section IV, the new optimization problem is formulated and solved for each CAV to address passenger discomfort inside the MZ. Furthermore, we investigate the associated tradeoff between minimizing fuel consumption and a measure of passenger discomfort. Finally, we present simulation results in Section V, and concluding remarks in Section VI.


\section{The Model}
\label{sec:1a}
\begin{figure}[ptb]
\centering
\includegraphics[width=3in]{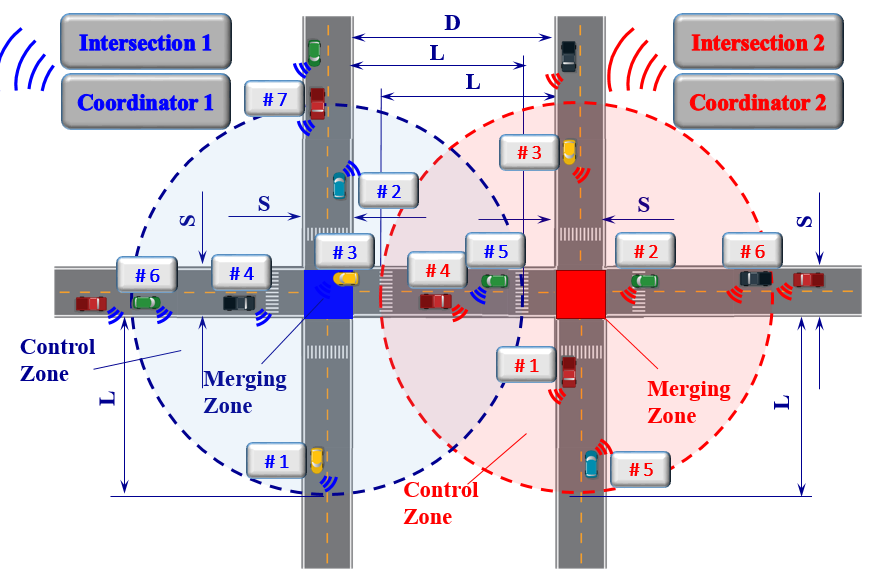} \caption{Connected and automated crossing two adjacent intersections.}%
\label{fig:intersection}%
\end{figure}

We briefly review the model introduced in
\cite{ZhangMalikopoulosCassandras2016} and \cite{Malikopoulos2016} where there
are two intersections, 1 and 2, located within a distance $D$ (Fig.
\ref{fig:intersection}). The region at the center of each intersection, called
\emph{Merging Zone} (MZ), is the area of potential lateral CAV collision.
Although it is not restrictive, this is taken to be a square of side $S$. Each
intersection has a \emph{Control Zone} (CZ) and a coordinator that can
communicate with the CAVs traveling within it. The distance between the entry
of the CZ and the entry of the MZ is $L>S$, and it is
assumed to be the same for all entry points to a given CZ.

Let $N_{z}(t)\in\mathbb{N}$ be the cumulative number of CAVs which have entered the CZ and formed a queue 
by time $t$, $z = 1,2$. The way the queue is formed is not restrictive to our analysis in the rest of the paper. When a CAV reaches the CZ of intersection $z$, the coordinator
assigns it an integer value $i=N_{z}(t)+1$. If two or more
CAVs enter a CZ at the same time, then the corresponding coordinator
selects randomly the first one to be assigned the value $N_{z}(t)+1$. 

For simplicity, we assume that
each CAV is governed by second order dynamics%
\begin{equation}
\dot{p}_{i}=v_{i}(t)\text{, }~p_{i}(t_{i}^{0})=0\text{; }~\dot{v}_{i}%
=u_{i}(t)\text{, }v_{i}(t_{i}^{0})\text{ given}\label{eq:model2}%
\end{equation}
where $p_{i}(t)\in\mathcal{P}_{i}$, $v_{i}(t)\in\mathcal{V}_{i}$, and
$u_{i}(t)\in\mathcal{U}_{i}$ denote the position of the CAV $i$ starting from the entry of the CZ, speed and
acceleration/deceleration (control input) of each CAV $i$. The sets $\mathcal{P}_{i}$,
$\mathcal{V}_{i}$ and $\mathcal{U}_{i}$ are complete and totally bounded sets of $\mathbb{R}$. These dynamics are
in force over an interval $[t_{i}^{0},t_{i}^{f}]$, where $t_{i}^{0}$ and $t_{i}^{f}$ are the
times that the vehicle $i$ enters the CZ  and exits the MZ of intersection $z$ respectively.

To ensure that the control input and vehicle speed are within a given
admissible range, the following constraints are imposed:
\begin{equation}%
\begin{split}
u_{i,min}  &  \leq u_{i}(t) \leq u_{i,max},\quad\text{and}\\
0  &  \leq v_{min} \leq v_{i}(t) \leq v_{max},\quad\forall
t\in\lbrack t_{i}^{0},t_{i}^{m}],
\end{split}
\label{speed_accel constraints}%
\end{equation}
where $t_{i}^{m}$ is the time that the vehicle $i$ enters the MZ. 
To ensure the absence of any rear-end collision throughout the CZ, we impose
the \emph{rear-end safety} constraint
\begin{equation}
s_{i}(t)=p_{k}(t)-p_{i}(t) \geq \delta,\quad\forall t\in\lbrack t_{i}%
^{0},t_{i}^{m}] \label{rearend}%
\end{equation}
where $\delta$ is the \emph{minimal safe distance} allowable and $k$ is the CAV physically ahead of $i$.

The objective of each CAV is to derive an optimal acceleration/deceleration, in
terms of fuel consumption, inside the CZ, i.e., over the time interval $[t_{i}^{0},t_{i}^{m}]$. In addition, we impose safety
constraints to avoid both rear-end and lateral collisions inside
the MZ. The conditions under which the rear-end collision avoidance constraint does not become active inside the CZ are provided in \cite{Zhang2016}.


\section{Vehicle Coordination and Control}


\subsection{Modeling Left and Right Turns}
\label{sec:2a}
In order to include left and right turns, we impose the following assumptions:

\begin{assumption}
\label{ass:sensor} Each vehicle $i$ has proximity sensors and can observe
and/or estimate local information that can be shared with other vehicles.
\end{assumption}

\begin{assumption}
\label{ass:dir} The decision of each vehicle $i$ on whether a turn is to be made at the MZ is known upon its entry in the CZ.
\end{assumption}

Let $d_i$ denote the decision of vehicle $i$ on whether a turn is to be made at the MZ, where $d_i = 0$ indicates left turn, $d_i = 1$ indicates going straight and $d_i = 2$ indicates right turn. 


%

Left and right turns need special attention in the context of safety while ensuring passenger comfort. We impose the following three \textit{Maximal Allowable Speed} limits inside the MZ: (1) $v_L^a$ for CAVs planning to make left turns, (2) $v_R^a$ for CAVs making right turns, and (3) $v^a$ for CAVs going straight.

The objective for each CAV is to make a safe turn while minimizing a measure of passenger discomfort. Discomfort due to a turn originates in the tangential and the centripetal forces. As the latter is perpendicular to the former, it cannot change the speed of the vehicle and it is the magnitude of the tangential force which must be maintained fixed so that passengers are minimally affected by the turn. Thus, we use \textit{jerk}, i.e., the rate of change of acceleration which results in speed vibrations, as a measure of passenger discomfort. 

%


\subsection{Vehicle Communication Structure}

\label{sec:2b}

\begin{figure*}[htb!]
\centering
\includegraphics[width= 2\columnwidth]{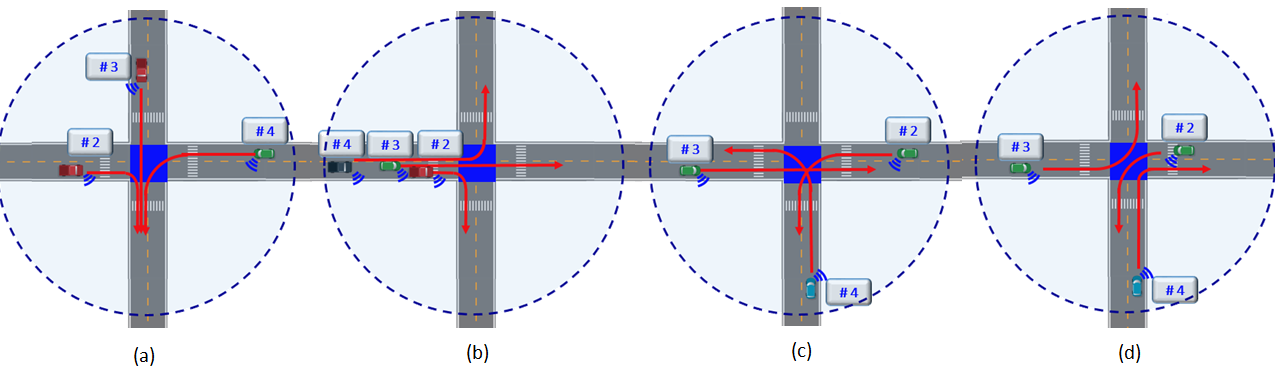} \caption{Illustration of different subsets of $\mathcal{Q}^{z}(t)$: (a) subset $\mathcal{E}^{z}(t)$; (b) subset $\mathcal{S}^{z}(t)$; (c) subset $\mathcal{L}^{z}(t)$; (d) subset $\mathcal{O}^{z}(t)$.}%
\label{category}%
\end{figure*}

When a CAV enters a CZ, $z=1,2$, it is assigned a set $\mathcal{Q}^{z}_{i}(t)$ from the coordinator, where  $\mathcal{Q}^{z}_{i}(t)$ = $\{\mathcal{E}^{z}_{i}(t),$ $\mathcal{S}^{z}_{i}(t),$ $\mathcal{L}^{z}_{i}(t),$ $\mathcal{O}^{z}_{i}(t)\}$, defined next, indicates the positional relationship between CAV $i$ and all other CAVs $j,$ $0<j<i$. With respect to CAV $i$, CAV $j$, $0<j<i$ belongs to one and only one of these subsets defined as follows:

(i) $\mathcal{E}^{z}_{i}(t)$ contains all vehicles that can cause rear-end collision at the \textit{end} of the MZ with $i$, e.g.,  $\mathcal{E}_{3}^{z}(t)$ contains vehicle \#2 as it may cause rear-end collision with vehicle \#3 at the end of the MZ (Fig. \ref{category}(a)), and $\mathcal{E}_{4}^{z}(t)$ contains vehicles \#3 and \#2 as it may cause rear-end collision with vehicle \#4 at the end of the MZ (Fig. \ref{category}(a)). Note that this subset does not contain the indices corresponding to vehicles cruising on the same lane and towards the same direction.

To ensure the absence of rear-end collision, the following constraint is applied:
\begin{equation}
t_{i}^{f}>t_{j}^{f},~ j \in \mathcal{E}^{z}_{i}(t).
\label{eq:fifo1}%
\end{equation}

(ii) $\mathcal{S}^{z}_{i}(t)$ contains all vehicles traveling on the same lane that can cause rear-end collision at the \textit{beginning} of the MZ with $i$, e.g., $\mathcal{S}_{3}^{z}(t)$ contains vehicle \#2 as it may cause rear-end collision with vehicle \#3 at the beginning of the MZ (Fig. \ref{category}(b)), and  $\mathcal{S}_{4}^{z}(t)$ contains vehicles \#3 and \#2 as it may cause rear-end collision with vehicle \#4 at the beginning of the MZ (Fig. \ref{category}(b)). Note that this subset contains the indices corresponding to vehicles cruising on the same lane and towards the same direction.

To ensure the absence of rear-end collision at the beginning of the MZ, the following condition is applied:
\begin{equation}
t_{i}^{m} > t_{k}^{m},~ k \in \mathcal{S}^{z}_{i}(t).
\label{eq:fifotim}%
\end{equation}

(iii) $\mathcal{L}^{z}_{i}(t)$ contains all vehicles traveling on different lanes and towards different lanes that can cause lateral collision inside the MZ with $i$, e.g.,  $\mathcal{L}_{3}^{z}(t)$ contains vehicle \#2 as it may cause lateral collision with vehicle \#3 inside the MZ (Fig. \ref{category}(c)), and  $\mathcal{L}_{4}^{z}(t)$ contains vehicles \#3 and \#2 as it may cause lateral collision with vehicle \#4 inside the MZ (Fig. \ref{category}(c)).

(iv) $\mathcal{O}^{z}(t)$ contains all vehicles traveling on different lanes and towards different directions that cannot cause lateral collision at the MZ with $i$, e.g., $\mathcal{O}_{3}^{z}(t)$ contains vehicle \#2 since it cannot cause any collision with vehicle \#3 (Fig. \ref{category}(d)), and $\mathcal{O}_{4}^{z}(t)$ contains vehicles \#3 and \#2 since it cannot cause any collision with vehicle \#4 (Fig. \ref{category}(d)).

Note that $\mathcal{Q}^{z}_{i}(t)$ defined above is different from the definition in \cite{ZhangMalikopoulosCassandras2016} where we had no turns. With left and right turns being considered, the collision scenarios are more complicated. As the positional relationship between CAV $i$ and $j$, $j < i$, $j \neq i-1$ cannot be explicitly determined through CAV $i-1$ with turns being involved, the terminal conditions of CAV $i$ cannot solely depend on CAV $i-1$.

Recalling that $t_i^m$ is the assigned time for CAV $i$ to enter the MZ, we require the following condition:
\begin{equation}
t_{i}^{f}\geq t_{i-1}^{f},~i>1.
\label{eq:fifo}%
\end{equation}
Note that in \cite{ZhangMalikopoulosCassandras2016} the speed in the MZ was considered to be constant, hence  $t_i^m > t_{k}^m$ could be ensured by $t_i^f > t_k^f$. However, with turns being considered, the speed and the trajectories in the MZ may be different for CAV $i-1$ and $i$, therefore, $t_i^f \geq t_{i-1}^f$ does not imply $t_i^m \geq t_{i-1}^m$. This becomes an issue only when $i-1 \in  \mathcal{S}^{z}_{i}(t)$. In that case, we have to ensure \eqref{eq:fifotim} holds. Also, for $i-1 \in  \mathcal{E}^{z}_{i}(t)$, \eqref{eq:fifo1} must be satisfied.

There is a number of ways to satisfy \eqref{eq:fifo}.  For example, we
may impose a strict First-In-First-Out (FIFO) queueing structure,
where each vehicle must enter the MZ in the same order
it entered the CZ. The crossing sequence may also be determined in a priority-based fashion. More generally,  $t_i^f$ (and $t_i^m$) may be determined for each vehicle $i$ at time $t_{i}%
^{0}$ when the vehicle enters the CZ. If \eqref{eq:fifo} is satisfied and \eqref{eq:fifo1} and \eqref{eq:fifotim} both hold, then the order in the queue is preserved. Otherwise, \eqref{eq:fifo} is violated, then the order may need to be updated so that CAV $i$ is placed
in the $j$th, $j<i$, queue position such that  \eqref{eq:fifo} is satisfied, and both \eqref{eq:fifo1} and \eqref{eq:fifotim} still hold. The policy through which the order
(\textquotedblleft schedule\textquotedblright) is specified may be the result
of a higher level optimization problem as long as the condition \eqref{eq:fifo}-\eqref{eq:fifotim} are preserved. In what
follows, we will adopt a specific scheme for determining $t_{i}^{m}$ and $t_i^f$ (upon
arrival of CAV $i$) based on our problem formulation, without affecting
the terminal conditions of 1, $\cdots$,$i-1$, but we emphasize that our analysis is not
restricted by the policy designating the order of the vehicles within the
queue.

For each CAV $i$, we define its
\emph{information set} $Y_{i}(t)$, $t\in\lbrack t_{i}^{0},t_{i}^{f}]$, as
\begin{gather}
Y_{i}(t)\triangleq\Big\{p_{i}(t),v_{i}(t),\mathcal{Q}^{z}_{i}(t) ,z=1,2, s_{i}(t),t_{i}^{m}, t_i^f, d_i \Big\},\label{InfoSet}%
\end{gather}
where $p_{i}(t),v_{i}(t)$ are the traveling distance and speed of CAV $i$ inside the
CZ it belongs to, and $\mathcal{Q}^{z}_{i}(t) = \{\mathcal{E}_{i}%
^{z}(t),$ $\mathcal{S}_{i}^{z}(t),$ $\mathcal{L}_{i}^{z}(t),$ $\mathcal{O}%
_{i}^{z}(t)\}$, $z=1,2.$ indicates the positional relationship with respect to CAV $j$, $j < i$. The fourth element in $Y_{i}(t)$ is $s_{i}(t)=p_{k}(t)-p_{i}(t)$,
the distance between CAV $i$ and CAV $k$ which is immediately
ahead of $i$ in the same lane (the index $k$ is made available to $i$ by the
coordinator under Assumption \ref{ass:sensor}. $t_{i}^{m}$ and $t_i^f$ are the times targeted for
CAV $i$ to enter and exit the MZ respectively, whose evaluation is discussed next. The last element $d_i$, indicates whether $i$ is making a left or right turn, or going straight at the MZ, which becomes known once the vehicle enters the CZ (Assumption \ref{ass:dir}). Note that once CAV $i$ enters the CZ, then all information in $Y_{i}(t)$ becomes
available to $i$.

For safety, comfort and fuel efficiency, it is appropriate for vehicles to make turns at an intersection at low speeds. The speed for which an intersection curve is designed depends on speed limit, the type of intersection, and the traffic volume \cite{aashto2001policy}. Generally, the ``desirable time" $\Delta_i$ that a vehicle needs to make a turn at an intersection \cite{aashto2001policy} is 

\begin{equation}
\Delta_i =  
\left\{
\begin{array}
[c]{ll}%
\frac{R_i}{\sqrt{15R_i(0.01E + F)}}, & \mbox{if $d_i = 0, 2$},\\
\frac{S}{v^a} , & \mbox{if $d_i = 1$}, 
\end{array}
\right.
\label{deltat1}
\end{equation}
where $R_i$ is the centerline turning radius, 
 $E$ is the super-elevation, which is zero in urban conditions, and $F$ is the side friction factor. Therefore, the time $t_i^m$ that the vehicle $i$ enters the MZ is directly related to the time $t_i^f$ that the vehicle exits the MZ through $\Delta_i$:
\begin{equation}
t_i^f = t_i^m + \Delta_i. \label{deltat2}
\end{equation}
Note that $\Delta_i$ is different for left and right turns since the corresponding turning radii $R_i$ are different.

\subsection{Terminal Conditions}
\label{tim}
We now turn our attention to the terminal conditions, i.e., CAV $i$'s time, speed, and position for entering/exiting the MZ, of each vehicle $i$. 

(i) Let $e =\max \limits_{j} \{j \in \mathcal{E}^z_{i}(t)\}$.  In this case, CAV $e$ is immediately ahead of CAV $i$ in the FIFO queue that may cause rear-end collision at the end of the MZ. To avoid such rear-end collision, $e$ and $i$ should maintain a minimal safety distance $\delta$, by the time vehicle $i$ exits the MZ. For simplicity, we make the following assumption:

\begin{assumption}
\label{ass:cruise} For each vehicle $i$, the speed remains constant after the MZ exit for at least a length $\delta$.
\end{assumption}
Note that this assumption is simply made for math calculation purpose. 

Given the assumption above, we set
\begin{equation}
t_i^f = t_{e}^f + \frac{\delta}{v_{e}^f} \label{Case E}%
\end{equation}
where $t_{e}^f$ and $t_{i}^f$ is the time that vehicle $e$ and $i$ exits the MZ, and $v_{e}^f$ is the speed of the vehicle $e$ at the exit of the MZ. The terminal speeds $t_i^m$ and $t_i^f$ are set as follows:

\begin{equation}
v_{i}^f= v_{i}^m =\left\{
\begin{array}
[c]{ll}%
v_L^a , & \mbox{if $d_{i} = 0$},\\
v^a , & \mbox{if $d_{i} = 1$}, \\
v_R^a , & \mbox{if $d_{i} = 2$}, 
\end{array}
\right.
\label{def:vif}
\end{equation}
where $v_{i}^m$ is the speed of the vehicle $i$ at the entry of the MZ. Note that $v_e^f$ can also be determined through \eqref{def:vif}. In our earlier work \cite{ZhangMalikopoulosCassandras2016}, we considered a constant speed for the CAVs inside the MZ. To consider left and right turns, the speed is no longer constant. Therefore, we formulate a new optimization problem to address a measure of passenger discomfort within the MZ. For this problem, \eqref{def:vif} are the initial and terminal speeds for each CAV, and  $t_i^m$ is the initial time which can be evaluated according to \eqref{deltat2} given $t_i^f$. Note that $v_{i}^m$ in \eqref{def:vif} and $t_i^m$ are the terminal conditions for the decentralized optimal control problem in the CZ. 

(ii) Let $s =\max \limits_{j} \{j \in \mathcal{S}^z_{i}(t)\}$. In this case, CAV $s$ is immediately ahead of CAV $i$ in the FIFO queue that may cause rear-end collision at the beginning of the MZ. To guarantee the rear-end collision constraint does not become active we set, 

\begin{equation}
t_i^m =t_{s}^m+ \Delta_{s}^{\delta},
\end{equation}
where
\begin{equation}
\Delta_{s}^{\delta} =  
\left\{
\begin{array}
[c]{ll}%
\frac{\delta}{\sqrt{15R_i(0.01E + F)}}, & \mbox{if $d_{s} = 0, 2$},\\
\frac{\delta}{v^a} , & \mbox{if $d_{s} = 1$}
\end{array}
\right.
\label{delta1}
\end{equation}
is the time vehicle $s$ needs to travel a distance $\delta$ inside the MZ. The time $t_{i}^f$ that the vehicle $i$ will be exiting the MZ can be evaluated from \eqref{deltat2} while $v_{i}^f$ is defined in \eqref{def:vif}.

Let $S_L$ and $S_R$ denote the length of the left and right turn trajectories, respectively. If vehicle $s$ makes a left turn and vehicle $i$ makes a right turn, $S_L > S_R$, implies that  $t_i^f < t_{s}^f$, thus \eqref{eq:fifo} does not hold. In that case, we set $t_i^f = t_{s}^f$, and  $t_i^m$ is evaluated  according to \eqref{deltat2}. The re-evaluation of $t_i^m$ can only make $t_i^m$ larger, thus, \eqref{eq:fifotim} still holds. Hence, we have

\begin{equation}
t_i^f = \max\{t_{s}^m+ \Delta_{s}^{\delta}+\Delta_i,t_s^f\}
\label{Case S}
\end{equation}

(iii) Let $l =\max \limits_{j} \{j \in \mathcal{L}^z_{i}(t)\}$. In this case, CAV $l$ is  immediately ahead of CAV $i$ in the FIFO queue that could cause lateral collision with $i$ inside the MZ. We constrain the MZ to contain only $l$ or $i$ so as to avoid lateral collision. 

\begin{equation}
t_i^m = t_{l}^f.
\label{lem11}
\end{equation}
The time $t_i^f$ can be evaluated through \eqref{deltat2} and $v_i^f$ is defined in \eqref{def:vif}.

(iv) Let $o =\max \limits_{j} \{j \in \mathcal{O}^z_{i}(t)\}$. In this case, CAV $o$ is  immediately ahead of CAV $i$ in the FIFO queue that will not generate any collision with $i$ in the MZ, so we set
\begin{equation}
t_i^f = t_{o}^f.
\label{lem12}
\end{equation}
The time $t_i^m$ can be evaluated through \eqref{deltat2} and $v_i^f$ is defined in \eqref{def:vif}. 

In order to ensure the absence of any collision type, we set $t_i^f$ as follows:

\begin{equation}
t_i^f = \max \{t_{e}^f + \frac{\delta}{v_{e}^f},t_{s}^m+ \Delta_{s}^{\delta}+\Delta_i,t_s^f, t_{l}^f + \Delta_i, t_o^f \} 
\label{tif}
\end{equation}

Recall that in \cite{ZhangMalikopoulosCassandras2016}, $t_i^f$ and $t_i^m$ can be recursively determined through CAV $i-1$ and $k$. However, with left and right turns being considered, the positional relationship between $i$ and $j$, $j <i, j\neq i-1$ becomes more complicated. CAV $i$ now depends on four CAVs  $e$, $s$, $l$, and $o$. However, the essence of the recursive structure stays the same. It follows from (\ref{deltat2}) and (\ref{tif}) that $t_{i}^{f}$ and $t_i^m$ can
always be recursively determined from CAVs $e$, $s$, $l$, and $o$, which preserves simplicity in the solution and enables decentralization.

Although (\ref{deltat2}) through (\ref{tif}) provide a
simple recursive structure for determining $t_{i}^{m}$, the presence of the
control and state constraints \eqref{speed_accel constraints} may prevent
these values from being admissible. This may happen by
\eqref{speed_accel constraints} becoming active at some internal point during
an optimal trajectory (see \cite{Malikopoulos2016} for details). In addition,
however, there is a global lower bound to $t_{i}^{m}$, which depends on $t_{i}^{0}$ and on whether CAV $i$
can reach $v_{max}$ prior to $t_{i-1}^{m}$ or not: $(i)$ If CAV $i$ enters the
CZ at $t_{i}^{0}$, accelerates with $u_{i,max}$ until it reaches $v_{max}$ and
then cruises at this speed until it leaves the MZ at time $t_{i}^{1}$, it was
shown in \cite{ZhangMalikopoulosCassandras2016} that
\begin{equation}
t_{i}^{1}=t_{i}^{0}+\frac{L}{v_{max}}+\frac{(v_{max}-v_{i}^{0})^{2}%
}{2u_{i,max}v_{max}}.
\end{equation}
$(ii)$ If CAV $i$ accelerates with $u_{i,max}$ but reaches the MZ at $t_{i}^{m}$
with speed $v_{i}^{m}<v_{max}$, it was shown in
\cite{ZhangMalikopoulosCassandras2016} that
\begin{equation}
t_{i}^{2}=t_{i}^{0}+\frac{v_{i}(t_{i}^{m})-v_{i}^{0}}{u_{i,max}},
\end{equation}
where $v_{i}(t_{i}^{m})=\sqrt{2L u_{i,max}+(v_{i}^{0})^{2}}$. Thus,
\begin{equation}
t_{i}^{c}=t_{i}^{1} \mathds{1}_{ v_i^m = v_{max}} + t_{i}^{2} (1 - \mathds{1}_{ v_i^m = v_{max}} ) \nonumber
\end{equation}
is a lower bound of $t_{i}^{m}$ regardless
of the solution of the problem.

\subsection{Decentralized Control Problem Formulation and Analytical Solution}
\label{sec:2c}
Recall that at time $t$, the values of $t_{i-1}^{f}$, $t_{i-1}^m$, $v_{i-1}^f$,
$\mathcal{Q}^{z}_{i}(t)$, $z=1,2$, $d_{i-1}$ are available to CAV $i$ through its information set in (\ref{InfoSet}). This is
necessary for $i$ to compute $t_{i}^{f}$ and $t_i^m$ appropriately and satisfy
\eqref{eq:fifo} and (\ref{eq:fifotim}). 


The decentralized optimal control problem for each CAV approaching either
intersection is formulated so as to minimize the $L^{2}$-norm of its control
input (acceleration/deceleration). There is a monotonic relationship between fuel consumption for each
CAV $i$, and its control input $u_{i}$  \cite{Rios-Torres2015}. Therefore, we formulate the
following problem for each $i$:
\begin{gather}
\min_{u_{i}\in U_{i}}\frac{1}{2}\int_{t_{i}^{0}}^{t_{i}^{m}}K_{i}\cdot
u_{i}^{2}~dt\nonumber\\
\text{subject to}:\eqref{eq:model2},(\ref{speed_accel constraints}%
),t_i^m,\text{ }p_{i}(t_{i}^{0})=0\text{, }%
p_{i}(t_{i}^{m})=L,\label{eq:decentral}\\
z=1,2, \text{and given }t_{i}^{0}\text{, }v_{i}(t_{i}^{0}),\nonumber
\end{gather}
where $K_{i}$ is a factor to capture CAV diversity (for simplicity we set
$K_{i}=1$ for the rest of this paper). Note that this formulation does not include the safety constraint (\ref{rearend}).

An analytical solution of problem \eqref{eq:decentral} may be obtained through
a Hamiltonian analysis. The presence of constraints (\ref{speed_accel constraints}) and
\eqref{rearend} complicates this analysis. The complete solution including the constraints in (\ref{speed_accel constraints}) is given in
\cite{Malikopoulos2016}. Assuming that the constraints are
not active upon entering the CZ and that they remain inactive throughout
$[t_{i}^{0},t_{i}^{m}]$, a complete solution was derived in
\cite{Rios-Torres2015} and \cite{Rios-Torres2} for highway on-ramps, and in
\cite{ZhangMalikopoulosCassandras2016} for two adjacent intersections. The
solution to \eqref{eq:decentral} differs from what we have derived in \cite{ZhangMalikopoulosCassandras2016} in the values of the coefficients instead of the structure, as the terminal conditions from (\ref{deltat2}) through (\ref{lem12}) consider left and right turns. 
The optimal control input
(acceleration/deceleration) over $t\in\lbrack t_{i}^{0},t_{i}^{m}]$ is given
by
\begin{equation}
u_{i}^{\ast}(t)=a_{i}t+b_{i}\label{eq:20}%
\end{equation}
where $a_{i}$ and $b_{i}$ are constants. Using (\ref{eq:20}) in the CAV
dynamics \eqref{eq:model2} we also obtain the optimal speed and position:
\begin{equation}
v_{i}^{\ast}(t)=\frac{1}{2}a_{i}t^{2}+b_{i}t+c_{i}\label{eq:21}%
\end{equation}%
\begin{equation}
p_{i}^{\ast}(t)=\frac{1}{6}a_{i}t^{3}+\frac{1}{2}b_{i}t^{2}+c_{i}%
t+d_{i},\label{eq:22}%
\end{equation}
where $c_{i}$ and $d_{i}$ are constants of integration. The constants $a_{i}$,
$b_{i}$, $c_{i}$, $d_{i}$ can be computed by using the given initial and final
conditions. The analytical solution (\ref{eq:20}) is only
valid as long as all initial conditions satisfy (\ref{speed_accel constraints}) and
\eqref{rearend} and none of these constraints becomes active in
$[t_{i}^{0},t_{i}^{m}]$. Otherwise, the solution needs to be modified as described in \cite{Malikopoulos2016}.
Recall that the constraint (\ref{rearend}) is not
included in (\ref{eq:decentral}) and it is a much more challenging matter. To
address this, we derive the conditions under
which the CAV's state maintains feasibility in terms of satisfying \eqref{rearend} over $[t_{i}^{0},t_{i}^{m}]$ in \cite{Zhang2016}.

\section{Joint Minimization of Passenger Discomfort and Fuel Consumption in the MZ}

\subsection{Passenger Discomfort}
\label{sec:3a}

It is reported in \cite{hogan1984adaptive} that the comfort of the passengers in transportation can be quantified as a function of jerk, which is the time derivative of acceleration, i.e., $J_i(t) = \dot{u}_i(t)$.  
Hence, the following optimization problem is formulated with the objective of minimizing the $L^2$-norm of jerk for each vehicle $i$, where the acceleration/deceleration $u_i(t)$ is the control input:
\begin{gather}
\min_{u_{i}}\frac{1}{2}\int_{t_{i}^{m}}^{t_{f}^{m}} J_{i}^{2}~dt  \label{eq:comfort}\\
\text{subject to}:  \eqref{eq:model2}, ~ J_i(t) = \dot{u}_i(t),  \nonumber \\
 u_i(t_i^m), u_i(t_i^f), v_i^m, v_i^f, p_i(t_i^m),  p_i(t_i^f), \text{given  } t_i^m, t_i^f, \nonumber 
\end{gather}

The analytical solution of problem \eqref{eq:comfort} has been obtained in \cite{ntousakis2016optimal} using Hamiltonian analysis and considering the jerk as the control input. However, here we control jerk indirectly through the acceleration/deceleration, so the analytical closed-form solution is
\begin{equation}
u_{i}^{\ast}(t)=\frac{1}{6}a_{i}t^{3}+\frac{1}{2}b_{i}t^{2}+c_{i}t+d_{i},\label{22}%
\end{equation}
\begin{equation}
v_{i}^{\ast}(t)=\frac{1}{24}a_{i}t^{4}+\frac{1}{6}b_{i}t^{3}+ \frac{1}{2}c_{i}t^2+d_{i} t + e_i,\label{23}%
\end{equation}
\begin{equation}
p_{i}^{\ast}(t)=\frac{1}{120}a_{i}t^{5}+\frac{1}{24}b_{i}t^{4}+ \frac{1}{6}c_{i}t^3+ \frac{1}{2}d_{i} t^2 + e_i t + f_i,\label{24}%
\end{equation}
where $a_i$, $b_i$, $c_i$, $d_i$, $e_i$ and $f_{i}$ are constants of integration, which can be computed by using the given initial and final
conditions at $t_i^m$ and $t_i^f$.
Therefore, \eqref{22} is the analytical optimal control input corresponding to \eqref{eq:comfort} that will yield the minimum $L^2$-norm of jerk for vehicle $i$ inside the MZ.

\subsection{Tradeoff between Fuel Consumption and Passenger Discomfort}
\label{sec:3c}

To investigate this tradeoff between fuel consumption and passenger discomfort, we consider a convex combination of acceleration/deceleration and jerk to formulate the following optimization problem:

\begin{gather}
\min_{u_{i}}\frac{1}{2}\int_{t_{i}^{m}}^{t_{f}^{m}} (w\cdot q_1\cdot u_i^2 + (1- w)\cdot q_2\cdot J_{i}^{2})~dt  \label{eq:comfort1}\\
\text{subject to}:  \eqref{eq:model2}, ~J_i(t) = \dot{u}_i(t)  \nonumber \\
0 \leq w \leq 1, \nonumber \\
 u_i(t_i^m), u_i(t_i^f), v_i^m, v_i^f, p_i(t_i^m)\text{ and } p_i(t_i^f), \text{given  } t_i^m, t_i^f. \nonumber 
\end{gather}
where $q_1$, $q_2$ are normalization factors which are selected so that $q_1 \cdot u_i^2 \in [0,1]$ and $q_2 \cdot J_i^2 \in [0,1]$. The Hamiltonian function of \eqref{eq:comfort1} becomes
\begin{gather}
H_{i}\big(t, x(t),j(t)\big) = w \cdot q_1 \cdot \frac{1}{2} u_i^2 + (1 - w) \cdot q_2 \cdot \frac{1}{2} J^{2}_{i} \nonumber \\
+ \lambda^{p}_{i} \cdot
v_{i} + \lambda^{v}_{i} \cdot u_{i}  + \lambda^{u}_{i} \cdot J_{i}, \nonumber%
\end{gather}
where $\lambda^{p}_{i}$, $\lambda^{v}_{i}$ and $\lambda^{u}_{i}$ are the costates. Given the necessary conditions for optimality, we have the following second-order ordinary differential equation,
\begin{gather}
(1 - w) \cdot  q_2 \cdot  \ddot v -w \cdot  q_1 \cdot  v + \frac{1}{2}a_i t^2 + b_i t + c_i = 0 \nonumber
\end{gather}
from which, we can derive the optimal solution as

\begin{figure}[ptb]
\centering
\includegraphics[width= 1\columnwidth]{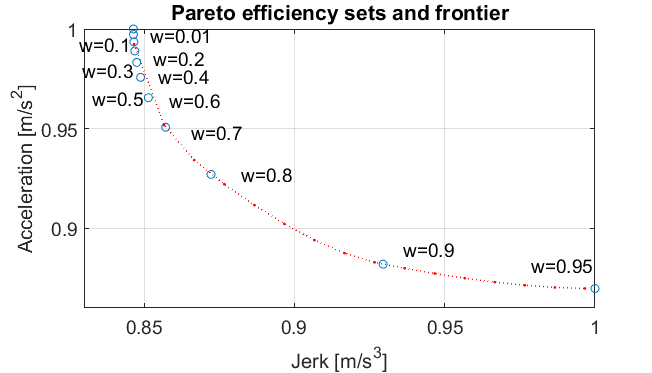} \caption{Pareto efficiency sets and frontier corresponding to different combinations of fuel consumption and passenger discomfort in the MZ.}%
\label{pareto}%
\end{figure}

\begin{gather}
J_i^*(t) = \frac{a_i}{w q_1} + e_i A_1^3 e^{A_1 t} + f_i A_2^3 e^{A_2 t}\\
u_i^*(t) = \frac{1}{wq_1}(a_i t + b_i) + e_i A_1^2 e^{A_1 t} + f_i A_2^2 e^{A_2 t}\\
v_i^*(t) = \frac{1}{wq_1} ( \frac{1}{2}a_i t^2 + b_i t + c_i  + \frac{a_i(1-w)q_2}{w q_1}) \nonumber\\
+ e_i A_1 e^{A_1 t} + f_i A_2 e^{A_2 t}\\
p_i^*(t) = \frac{1}{wq_1}(\frac{1}{6}a_i t^3 + \frac{1}{2}b_i t^2 + c_i t + \frac{a_i(1-w)q_2}{wq_1}t + d_i)  \nonumber  \\
+ e_i e^{A_1 t} + f_i e^{A_2 t}
\end{gather}

where

\begin{gather}
A_1 = \sqrt{\frac{w q_1}{(1 - w)q_2}}, ~~~~A_2 = -\sqrt{\frac{wq_1}{(1 - w)q_2}}. \nonumber 
\end{gather}
The constants
$a_i$, $b_i$, $c_i$, $d_i$, $e_i$ and $f_i$ can be computed using the initial and final conditions at $t_i^m$ and $t_i^f$. 
Note that since $ 0 \leq w \leq 1$, the optimal solution is only valid when $w \neq 1$ and $w \neq 0$. When $w = 0$, the problem reduces to \eqref{eq:comfort}. When $w = 1$, the problem becomes the same as the one formulated for the CZ in \eqref{eq:decentral}, which now minimizes the fuel consumption in the MZ:

\begin{gather}
\min_{u_{i}}\frac{1}{2}\int_{t_{i}^{m}}^{t_{i}^{f}}
u_{i}^{2}~dt \label{eq:fuel}\\
\text{subject to}:\eqref{eq:model2}, v_i^m, v_i^f, p_i(t_i^m),  p_i(t_i^f), \text{given  } t_i^m, t_i^f. \nonumber 
\end{gather}


To illustrate the tradeoff between passenger discomfort and fuel consumption, we examine a range of cases with different weights and produce the Pareto sets. The parameters used are listed in Sec. \ref{sim} except for the initial and terminal acceleration which are set to 0. By yielding all of the optimal solutions to \eqref{eq:comfort1} while varying the weight $w$, we can derive the Pareto sets and the Pareto frontier corresponding to different combinations of fuel consumption and passenger discomfort as shown in Fig. \ref{pareto}. Note that as $w\rightarrow0$, the solution of \eqref{eq:comfort1} becomes that of \eqref{eq:comfort}. 


\section{Simulation Examples}
\label{sim}

\begin{figure*}[htb!]
\centering
\includegraphics[width= 1.5\columnwidth]{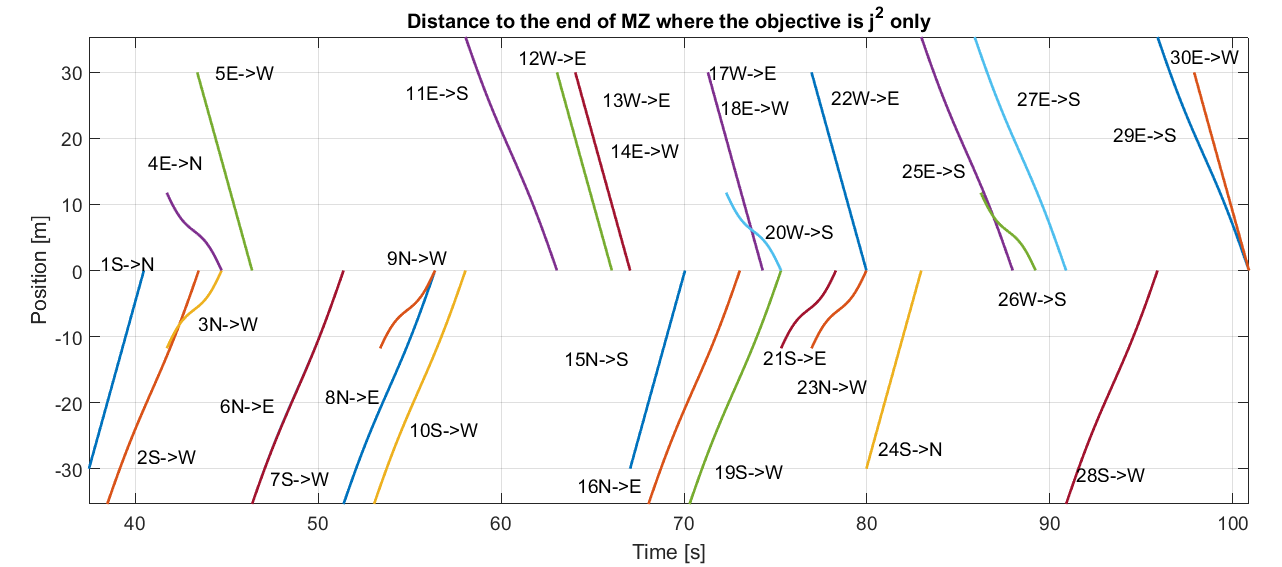} \caption{Distance to the end of MZ of the first 30 CAVs in the MZ.}%
\label{p}%
\end{figure*}

\begin{figure*}
\centering
\includegraphics[width= 1.7 \columnwidth]{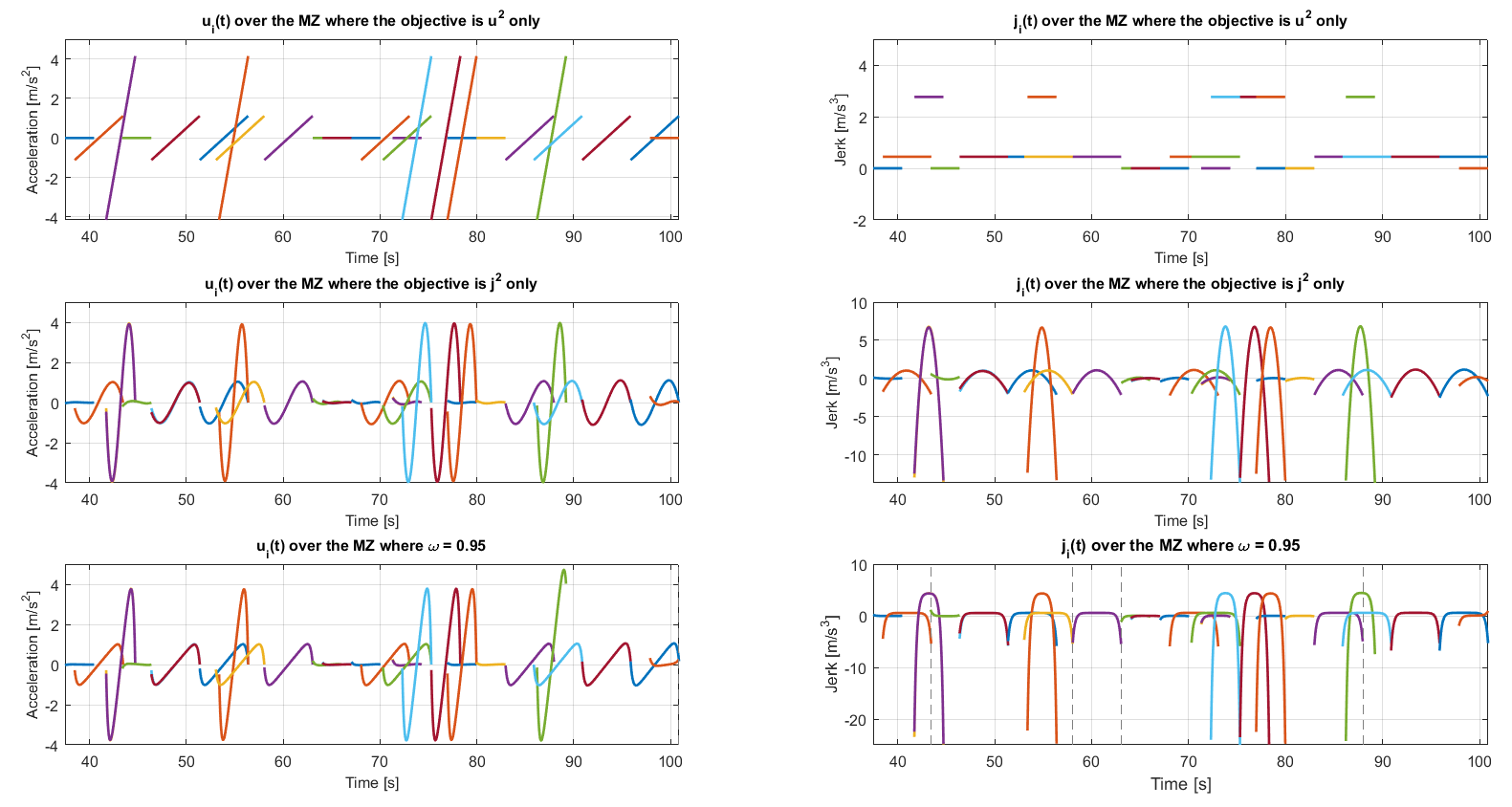} \caption{Acceleration/deceleration $u_i(t)$ and jerk $J_i(t)$ trajectories for the cases with different objectives: (a) minimize fuel consumption only; (b) minimize passenger discomfort only; (c) minimize a weighted sum of fuel consumption and passenger discomfort where $w = 0.95$.}%
\label{uj}%
\end{figure*}

The proposed decentralized optimal control framework incorporating turns is illustrated through simulation in MATLAB. For each direction, only one lane is considered. The parameters used are: $L=400$ m, $S=30$ m, $S_L = \frac{3}{8}\pi S$, $S_R = \frac{1}{8}\pi S$, $\delta=10$ m, $v_L^a = 8$ m/s, $v_R^a = 6$ m/s, $v^a = 10$ m/s and $\Delta_i = 5, 3, 3$s for left turn, going straight and right turn respectively. CAVs arrive at the CZ based on a random arrival process. Here, we assume a Poisson arrival process with rate $\lambda = 1$ and the speeds are uniformly distributed over $[10, 12]$.

We first consider the case where only the $L^2$-norm of jerk is optimized over the MZ. The initial and terminal conditions of time and speed are defined from (\ref{deltat2}) through (\ref{lem12}). Observe that $p_i(t_i^m) = L$ and $p_i(t_i^f) = L + S$ if $i$ is going straight ($p_i(t_i^f) = L + S_L$ for left turn and $p_i(t_i^f) = L + S_R$ for right turn). The
two additional conditions needed for acceleration/deceleration are set as follows: (a) the initial acceleration for \eqref{eq:comfort} is set to the terminal value derived from the optimal control problem in the CZ \eqref{eq:decentral}, under which the acceleration/deceleration is continuous at $t_i^m$; (b) the terminal acceleration is set as zero. The position trajectories of the first 30 CAVs in the MZ are shown in Fig. \ref{p}. CAVs are separated into two groups: CAVs shown above zero
are driving from east or west, and those below zero are
driving from north or south, with labels indicating the position of the vehicles in the FIFO queue and the driving direction. These figures include
different instances from each of Cases 1), 2), 3) or 4) in Sec. \ref{sec:2b}
regarding the value of $t_{i}^{f}$. For example, CAV \#11 is assigned
$t_{11}^{m}=t_{10}^{f}$, which corresponds to Case 3),
CAV \#23 is assigned $t_{23}^{m}=t_{22}^{m}$, which corresponds to Case 4), CAV \#27 is assigned $t_{27}^f = t_{26}^f + \frac{\delta}{26}$, which corresponds to Case 1), whereas CAV \#13 is assigned $t_{13}^m = t_{12}^m + \Delta_{12}^{\delta}$, which corresponds to Case 2).

To demonstrate the effectiveness of our optimal control in the MZ and compare formulations with different objectives, we examine three cases: (a) the objective is to minimize the $L^2$-norm of acceleration/deceleration only \eqref{eq:fuel}; (b) the objective is to minimize the $L^2$-norm of jerk only \eqref{eq:comfort}; (c) the objective is to minimize the weighted sum of $L^2$-norm of acceleration/deceleration and  $L^2$-norm of jerk \eqref{eq:comfort1}, where $w = 0.95$. Note that all terms should be normalized into a uniform, dimensionless scale for multi-objective optimization. The acceleration and jerk profiles of the first 30 CAVs are shown in Fig. \ref{uj}. Note that the optimal solution depends on how we set the initial and terminal acceleration/deceleration.

%
%
\section{Conclusions}
Earlier work \cite{ZhangMalikopoulosCassandras2016} has established a decentralized optimal control framework for optimally controlling CAVs crossing two adjacent intersections in an urban area. In this paper, we extended the solution of this problem  to account for left and right turns. In addition, we formulated and solved another optimization problem to minimize a measure of passenger discomfort while the vehicle turns at the intersection and investigated the associated tradeoff between minimizing fuel consumption and passenger discomfort. The optimal solution including turns do not require any additional computational time than what is required by the solution in \cite{ZhangMalikopoulosCassandras2016} since the terminal conditions are determined based on another set of collision-avoidance constraints, which can still enable online implementation. Future research should investigate the implications of having information with errors and/or delays to the system behavior.


\bibliographystyle{IEEETran}
\bibliography{CDC2017}

\end{document}